\documentclass{amsart}
\usepackage{amsthm,amscd,amssymb,verbatim,epsf,amsmath,eurosym,amsfonts,mathrsfs,graphicx,wrapfig}
\usepackage[colorlinks=true,linkcolor=blue,citecolor=blue]{hyperref}
\usepackage{amsmath, amsthm, amscd, amsfonts, amssymb, graphicx, color}
\makeatletter
\@namedef{subjclassname@2020}{\textup{2020} Mathematics Subject Classification}
\makeatother

\begin{document}

\makeatletter
\newcommand{\verbatimfont}[1]{\renewcommand{\verbatim@font}{\ttfamily#1}}
\makeatother
\newcommand{\R}{\mathbb{R}}
\newcommand{\C}{\mathbb{C}}
\newcommand{\om}{\omega}
\newcommand{\I}{{\mbox{I}}}
\newcommand{\II}{{\mbox{II}}}
\newcommand{\III}{{\mbox{III}}}
\newcommand{\Hess}{\mbox{\rm{Hess}}}
\newcommand{\Sing}{\mbox{\rm{Sing}}}
\newcommand{\Ind}{\mbox{\rm{Ind}}}
\theoremstyle{plain}
\newtheorem{Thm}{Theorem}
\newtheorem{Cor}{Corollary}
\newtheorem{Ex}{Example}
\newtheorem{Con}{Conjecture}
\newtheorem{Main}{Theorem 3}
\newtheorem{Lem}{Lemma}
\newtheorem{Prop}{Proposition}

\theoremstyle{definition}
\newtheorem{Def}{Definition}
\newtheorem{Note}{Note}
\newtheorem{Question}{Question}

\newtheorem{remark}{Remark}
\newtheorem{notation}{Notation}
\renewcommand{\thenotation}{}

\renewcommand{\rm}{\normalshape}%

\title[An index bound for smooth umbilic points]%
   { An index bound for smooth umbilic points}

\author{Brendan Guilfoyle}
\address{Brendan Guilfoyle\\
          School of STEM\\
          Munster Technological University\\
          Tralee \\
          Co. Kerry \\
          Ireland.}
\email{brendan.guilfoyle@mtu.ie}
\author{Wilhelm Klingenberg}
\address{Wilhelm Klingenberg\\
 Department of Mathematical Sciences\\
 University of Durham\\
 Durham DH1 3LE\\
 United Kingdom}
\email{wilhelm.klingenberg@durham.ac.uk }

\date{\today}

\begin{abstract}
We prove that the ${\mathbb Z}/2$-valued local index of an isolated umbilic point on a $C^{3+\alpha}$-smooth convex surface in Euclidean 3-space is less than two. The approach is to study the co-kernel of an associated Riemann-Hilbert boundary value problem. 

The link between the local and global is a semi-local technique that we term {\it totally real blow-up}.  Topologically, given a real surface in a complex surface, the totally real blow-up is the connect sum of the real surface with an embedded real projective plane.  We show that this increases the sum of the complex indices of the real surface by one, and hence cancels isolated hyperbolic complex points.

This leads to the reduction of the local result to the global result (the non-existence of closed embedded Lagrangian surfaces with a single complex point), proving that the umbilic index for smooth surfaces is less than two. Comparison of our smooth result with that of Hans Hamburger in the real analytic case (stating that the index of an isolated umbilic point on a real analytic convex surface is less than or equal to one) suggests the existence of ``exotic'' umbilic points of index $3/2$.
\end{abstract}

\maketitle
The purpose of this paper is to prove the following:

\vspace{0.1in}
\begin{Thm}\label{t:1}
    The index of any isolated umbilic point on a $C^{3+\alpha}$-smooth convex surface in Euclidean 3-space is less than two.
\end{Thm}
\vspace{0.1in}

Here the index (an element of ${\mathbb Z}/2$) is the winding number of the principal foliation about the isolated umbilic point, and $C^{3+\alpha}$ is the usual H\"older space with $\alpha\in(0,1)$. This bound does not preclude the existence of a smooth umbilic of index $3/2$, which is ruled out by the work of Hamburger in the real analytic case \cite{Ham}. We would argue that this points to a difference between the smooth and real analytic categories for surfaces in ${\mathbb R}^3$.

This resolves a well-known  Conjecture of Carath\'eodory \cite{Ham22}:

\begin{Thm}\label{t:2}
    The number of umbilic points on a $C^{3+\alpha}$-smooth convex 2-sphere in Euclidean 3-space is greater than one.
\end{Thm}

The prove starts with the reformulation of questions regarding umbilic points on convex surfaces in ${\mathbb R}^3$ in terms of questions regarding complex points on Lagrangian sections in $TS^2$ with its canonical neutral K\"ahler structure \cite{gak05}. Indeed, it is easily seen that Theorem \ref{t:1} is equivalent to:

\begin{Thm}\label{t:3}
    The index of any isolated complex point on a $C^{2+\alpha}$-smooth Lagrangian section $\Sigma$ of $TS^2$ is less than four.
\end{Thm}

The section $\Sigma$ is the set of oriented normal lines to the surface $S$, considered as a surface (with a loss of one derivative) in the 4-manifold $TS^2$ of all oriented lines in Euclidean 3-space. Theorem \ref{t:3} is proven as follows. 

Suppose for the sake of contradiction that there exists a $C^{2+\alpha}$-smooth Lagrangian section of $TS^2$ containing an isolated complex point of index $4+k$ for $k\geq 0$. Join a neighbourhood containing the complex point to a totally real Lagrangian hemisphere by the addition of a Lagrangian annulus. 

Generically, this annulus will contain isolated elliptic and hyperbolic complex points, where the complex index is $+1$ and $-1$, respectively. By Lai's index formula (see e.g. \cite{forst}) the sum of the complex indices on $\Sigma$ is $\chi(T\Sigma)+\chi(N\Sigma)=4$, and so the sum of the hyperbolic and elliptic indices is $-k$.

A straight-forward application of the h-principle proves that elliptic and hyperbolic complex points can be cancelled pair-wise in the Lagrangian category, leaving exactly $k$ hyperbolic points in the annulus. Denote the resulting Lagrangian section by $\Sigma'$. A corresponding surface in ${\mathbb R}^3$ that is orthogonal to the oriented lines of  $\Sigma'$ is shown schematically below with its umbilics and their indices.

\begin{wrapfigure}{l}{0.5\textwidth}
  \begin{center}
    \includegraphics[width=0.48\textwidth]{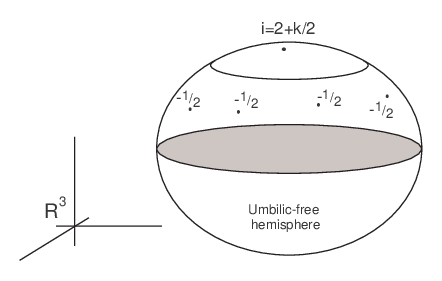}
  \end{center}
\end{wrapfigure}

We show that the connect sum of the Lagrangian section $\Sigma'$ in $TS^2$ with $k$ copies of ${\mathbb R}P^2$ removes hyperbolic complex points. Indeed, this construction can be applied to remove hyperbolic complex points on any real surface in a complex surface. The analogy with blowing-up in algebraic geometry motivates us to call it a {\it totally real blow-up}, where, in contrast to the standard notion of blowing up, the ambient complex surface is left unchanged by our operation.
 
Removal of the hyperbolic points yields $\Sigma_1=\Sigma'\#k{\mathbb R}P^2\subset TS^2$, a compact embedded surface containing a single complex point and a totally real Lagrangian hemisphere. Now by Theorem 68 of \cite{gak24} any closed totally real Lagrangian hemisphere admits a holomorphic disc with boundary lying upon it, implying that the co-kernel of the $\bar{\partial}$-operator is non-zero. 

On the other hand, the infinite dimensional Sard-Smale theorem for global sections of plane bundles with a single complex point implies that the co-kernel is zero \cite{gak19}. This result is extended to spheres with cross-caps in section 3.4 and the ensuing contradiction implies that a Lagrangian section containing a complex point of index $4+k$  cannot exist.

In the next section we summarize the background geometry required to prove the Theorem \ref{t:3}. In section \ref{s:2} we show how to remove hyperbolic complex points by connect sum with embedded real projective planes. Section \ref{s:3} contains the proof of the Theorem \ref{t:3} and the final section briefly discusses the result and its implications.

\section{Background} \label{s:1}

In this section we summarize the geometry of neutral K\"ahler $TS^2$. Further details can be found in \cite{gak05} and \cite{gak04}.

The space of oriented lines in ${\mathbb R}^3$ can be identified with  $TS^2$ by noting that an oriented line can be identified with a pair of orthogonal vectors $(\vec{U},\vec{V})$, the first of which is a unit vector (the direction of the oriented line):
\[
\{(\vec{U},\vec{V})\in{\mathbb R}^3\times{\mathbb R}^3\;|\;|\vec{U}|=1 \quad \vec{U}\cdot\vec{V}=0\;\}=TS^2.
\]
By lifting the standard complex coordinate $\xi$ (obtained by stereographic projection from the south pole on $S^2$) we get complex coordinates 
$(\xi,\eta)$ on $T^{10}S^2\cong TS^2$. In particular, identify $(\xi,\eta)\in{\mathbb {C}}^2$ with the vector
\[
\eta\frac{\partial}{\partial \xi}+\bar{\eta}\frac{\partial}{\partial \bar{\xi}}\in T_\xi S^2.
\]
In other words, the coordinate $\xi$ represents the direction  vector ${\vec{U}}$ of the oriented line, while $\eta$ determines the perpendicular distance vector ${\vec{V}}$. The canonical projection $\pi:TS^2\rightarrow S^2$, $\pi(\xi,\eta)=\xi$, maps an oriented line to its direction.

The neutral K\"ahler structure on $TS^2$ consists of a complex structure ${\mathbb J}$ for which these coordinates are holomorphic:
\[
{\mathbb J}\left(\frac{\partial}{\partial \xi}\right)=i\frac{\partial}{\partial \xi}
\qquad\qquad
{\mathbb J}\left(\frac{\partial}{\partial \eta}\right)=i\frac{\partial}{\partial \eta},
\]
together with a compatible symplectic 2-form $\Omega$ and a metric ${\mathbb G}$ of signature $(2,2)$, which have the following local expressions
in $(\xi,\eta)$-coordinates:
\[
\Omega=4(1+\xi\bar{\xi})^{-2}{\mathbb{R}}\mbox{e}\left(d\eta\wedge d\bar{\xi}-\frac{2\bar{\xi}\eta}{1+\xi\bar{\xi}}d\xi\wedge d\bar{\xi}\right),
\]
\[
{\mathbb{G}}=4(1+\xi\bar{\xi})^{-2}{\mathbb{I}}\mbox{m}\left(d\bar{\eta} d\xi+\frac{2\bar{\xi}\eta}{1+\xi\bar{\xi}}d\xi d\bar{\xi}\right).
\]
\vspace{0.1in}

Consider now a real surface $\Sigma$ in $TS^2$. A point $\gamma$ is said to be {\it complex} if the complex structure of $TS^2$ preserves the tangent space of $\Sigma$ at $\gamma$: ${\mathbb J}:T_\gamma\Sigma\rightarrow T_\gamma\Sigma$. In coordinates, a point $\gamma$ on $\Sigma$ given by $\nu\mapsto(\xi(\nu,\bar{\nu}),\eta(\nu,\bar{\nu}))$  is complex iff
\[
\partial\eta\bar{\partial}\xi-\bar{\partial}\eta\partial\xi|_{\gamma}=0.
\]
On the other hand a surface  $\Sigma$ is  {\it Lagrangian} at a point $\gamma$ if $\Omega|_\Sigma=0$ at $\gamma$. A Lagrangian surface is a surface which is Lagrangian at all of its points.

\begin{Prop}
    Let $S$ be a convex surface in ${\mathbb R}^3$ and $\Sigma\subset TS^2$ be the surface formed by the oriented normals to $S$. Then $\Sigma$ is a Lagrangian section of $\pi:TS^2\rightarrow S^2$.

Conversely, if $\Sigma$ is a Lagrangian section, then there exists a 1-parameter family of convex surfaces in ${\mathbb R}^3$ which are orthogonal to the oriented lines of $\Sigma$.

Moreover, $p\in S$ is umbilic iff the oriented normal to $S$ through $p$ is a complex point on the surface $\Sigma$ formed by the oriented normals of $S$.
\end{Prop}

Of particular importance are the surfaces in $TS^2$ that arise as a graph of a section of the projection $\pi:TS^2\rightarrow S^2$:

\begin{Prop}\label{p:support}
A section $\xi\mapsto(\xi,\eta=F(\xi,\bar{\xi}))$ is Lagrangian iff there exists a real function (called the {\em support function} of the corresponding surface in ${\mathbb R}^3$) $r:S^2\rightarrow {\mathbb R}$ such that
\[
\frac{\partial r}{\partial \xi}=\frac{2\bar{F}}{(1+\xi\bar{\xi})^2},
\]
where the support function is defined up to $r\rightarrow r+C$ (which yields parallel surfaces).

A point is complex iff at the point $\bar{\partial}F=0$.

A point on a Lagrangian section is complex iff the corresponding point on an orthogonal surface in ${\mathbb R}^3$ is umbilic.
\end{Prop}

An isolated umbilic point $p$ on a convex surface $S$ is a singularity of the principal foliation of the surface and, as such, has an index $i(p)\in{\mathbb Z}/2$ (a half-integer because the foliation may not be orientable).

On the other hand, an isolated complex point $\gamma$ on a real surface $\Sigma$ in a complex surface also has an index $I(\gamma)\in{\mathbb Z}$, see \cite{forst}. In our case, these are related by

\begin{Prop}\cite{gak04}
Let $S$ be a convex surface in ${\mathbb R}^3$ containing an isolated umbilic point $p$ and let $\Sigma\subset TS^2$ be the surface determined by the oriented normal lines of $S$ with corresponding isolated complex point $\gamma$. Then
\[
I(\gamma)=2i(p).
\]
\end{Prop} 

Following the standard convention for complex points we adopt the terminology:

\begin{Def}
An isolated complex point on a real surface is {\it elliptic} if the index is equal to one, and is {\it hyperbolic} if the index is equal to $-1$. 
\end{Def}

\vspace{0.1in}

\section{Totally Real Blow-up}\label{s:2}

Consider the connected sum of a real surface $\Sigma$ in $TS^2$ with a copy of ${\mathbb R}P^2$. That is, remove discs from $\Sigma$ and ${\mathbb R}P^2$ and identify the boundary circles. A copy of the real projective plane with a disc removed is called a {\it cross-cap}. It can also be viewed as an annulus with the inner boundary curve antipodally identified.

We term this operation {\it totally real blow-up}: ``blow-up''  because, at the topological level, it is the real analogue of complex blowing-up 
(connect sum of a complex surface with $\overline{{\mathbb C}P^2}$), and ``totally real'' because it removes certain types of complex points. 
Exactly which type can be removed is established in the next result:

\begin{Prop}\label{p:trb}
Let $\Sigma\subset{\mathbb M}$ be an embedded real surface in a complex surface with a single hyperbolic complex point $\gamma\in\Sigma$. Then the 
surface $\tilde{\Sigma}=\Sigma\#{\mathbb R}P^2$ given by removing a neighbourhood of $\gamma$ and joining a cross-cap can be smoothed so 
that $\tilde{\Sigma}$ is totally real and embedded. 
\end{Prop}
\begin{proof}
Choose local holomorphic coordinates $(\xi,\eta)$ on ${\mathbb M}$ so that the surface in the neighbourhood of the complex point at $\xi=0$  is given by
\[
\eta=\alpha\bar{\xi}^2+\beta\xi\bar{\xi}+o(|\xi|^3),
\]
for complex numbers $\alpha,\beta$. As the complex point is hyperbolic we have that $|\alpha|>2|\beta|$ (see \cite{forst}) and so by a 
compactly supported deformation this can be reduced to the form
\[
\eta=\alpha\bar{\xi}^2,
\]
without the creation of any further complex points.

Fix constants $\epsilon$ and  $R_0$ such that $1-\epsilon<R_0<1$, and consider the surface $\tilde{\Sigma}$ defined by
\begin{align}
\xi=&\;\;(1-\nu\bar{\nu})\nu \qquad\qquad {\mbox{ for}}\qquad 1-\epsilon\leq|\nu|\leq 1,\nonumber\\
& \nonumber\\
\eta=&\left\{
\begin{array}{lcl}
\alpha(1-\nu\bar{\nu})^2\bar{\nu}^2  &{\mbox{ for}} & 1-\epsilon\leq|\nu|\leq R_0 ,\\
 & \\
\alpha(a+b(1-\nu\bar{\nu})+c(1-\nu\bar{\nu})^2)\bar{\nu}^2&{\mbox{ for}}& R_0\leq|\nu|\leq 1.
\end{array}
\right.\nonumber
\end{align} 
This surface is $C^1$-smooth for the following choice of constants $a$ and $b$:
\[
a=(c-1)(1-R_0^2)^2 \qquad b=2(1-c)(1-R_0^2),
\]
and $c$ to be determined.  Moreover, it is the connected sum $\tilde{\Sigma}=\Sigma\#{\mathbb R}P^2$ and is easily seen to be embedded.

We now show that for a certain range of values of the constant $c$, the surface $\tilde{\Sigma}$ is totally real. Recall, a point $\gamma$ on a real surface is 
complex iff
\[
\partial\eta\bar{\partial}\xi-\bar{\partial}\eta\partial\xi|_{\gamma}=0.
\]
Computing this for the surface above, we find that for $ 1-\epsilon\leq|\nu|\leq R_0$ there are no complex points, while for $R_0\leq|\nu|\leq 1$ we have
\begin{align}
\partial\eta\bar{\partial}\xi-\bar{\partial}\eta\partial\xi=&-2[1+(c-1)R_0^4
-(c(3+2R_0^2)R_0^2+(5+2R_0^2)(1-R_0^2))\nu\bar{\nu}\nonumber\\
&\qquad\qquad\qquad\qquad\qquad+(5(1-R_0^2)+c(2+5R_0^2))\nu^2\bar{\nu}^2
  +3c\nu^3\bar{\nu}^3]\bar{\nu}.\nonumber
\end{align}
This is zero when $\nu=0$, which is outside of the range for the parameter, so we exclude this value.

Consider then the other factor with $x=\nu\bar{\nu}$ and $y=R_0^2$:
\[
g(x,y):=1+(c-1)y^2-(c(3+2y)y+(5+2y)(1-y))x+(5(1-y)+c(2+5y))x^2+3cx^3
\]
The function $g$ has the following properties:
\begin{enumerate}
\item $g(1,1)=0$,
\item $\partial_{x}g|_{(1,1)}=0$ and $\partial_{y}g|_{(1,1)}=0$,
\item $\partial^2_{x}g|_{(1,1)}=4c$ , $\partial^2_{y}g|_{(1,1)}=2(c-1)$ and $\partial_x\partial_{y}g|_{(1,1)}=-3(c-1)$.
\end{enumerate}
Thus $\partial^2_{x}g\partial^2_{y}g-(\partial_x\partial_{y}g)^2|_{(1,1)}=(9-c)(c-1)$ and so $(1,1)$ is local minimum for $g$ if $1<c<9$. 
Thus for this range of $c$, there exists $\epsilon>0$ s.t. $g(x,y)>0$ for $1-\epsilon< y^2<1$ and $R_0\leq x^2\leq 1$. For $c$ outside of this range, 
circles of complex points occur which contribute zero to the total complex index and may be removed by a generic perturbation.

We conclude that $\tilde{\Sigma}$ is totally real.
\end{proof}

To get a picture of the totally real blow-up of a hyperbolic point in terms of oriented lines in ${\mathbb R}^3$, consider the 
real surface in $TS^2$ given by the section $\xi\mapsto(\xi,\eta=(1+\xi\bar{\xi})^2\bar{\xi}^2)$. It is easily checked that this
surface is Lagrangian and has a single hyperbolic complex point at $\xi=0$:
\[
{\mathbb I}{\mbox m}\frac{\partial}{\partial \xi}\left(\frac{\eta}{(1+\xi\bar{\xi})^2}\right)=0
\quad
\frac{\partial \eta}{\partial \bar{\xi}}(0)=0
\quad
i(0)=\lim_{R\rightarrow 0}{\textstyle{\frac{1}{8\pi i}}}
      \int_0^{2\pi}\frac{\partial}{\partial \theta}\ln\left(\frac{\bar{\partial}\eta}{\partial\bar{\eta}}\right)d\theta=-{\textstyle{\frac{1}{2}}},
\]
where $\xi=Re^{i\theta}$. Thus the lines are perpendicular to a 1-parameter family of convex surfaces in ${\mathbb R}^3$ which contain
an isolated umbilic point of index $-1/2$. Let us construct these surfaces explicitly.

The support function is found by integrating the defining equation in Proposition \ref{p:support} and the result is
\[
r={\textstyle{\frac{2}{3}}}(\xi^3+\bar{\xi}^3)+C.
\]
To reconstruct these surfaces recall the fundamental correspondence relation \cite{gak04} in Euclidean coordinates $(x^1,x^2,x^3)$:
\[
x^1+ix^2=\frac{2(\eta-\bar{\eta}\xi^2)+2\xi(1+\xi\bar{\xi})r}{(1+\xi\bar{\xi})^2},
\qquad
x^3=\frac{-2(\eta\bar{\xi}+\bar{\eta}\xi)+(1-\xi^2\bar{\xi}^2)r}{(1+\xi\bar{\xi})^2},
\]
Therefore the 1-parameter family of convex surfaces, parameterized by the inverse of their Gauss maps, is
\[
x^1+ix^2=\frac{2(3\bar{\xi}^2-\xi^4+5\xi\bar{\xi}^3-3\xi^5\bar{\xi})}{3(1+\xi\bar{\xi})}+\frac{2c\xi}{1+\xi\bar{\xi}}
\]
\[
x^3=-\frac{4(\xi+\bar{\xi})(\xi-\xi\bar{\xi}+\xi)(1+2\xi\bar{\xi})}{3(1+\xi\bar{\xi})}+\frac{c(1-\xi\bar{\xi})}{1+\xi\bar{\xi}}
\]
Below is a plot of one of these surfaces with curves indicating the images of $|\xi|$=constant and of arg$(\xi)$=constant.

\begin{wrapfigure}{r}{0.5\textwidth}
  \begin{center}
\includegraphics[width=0.48\textwidth]{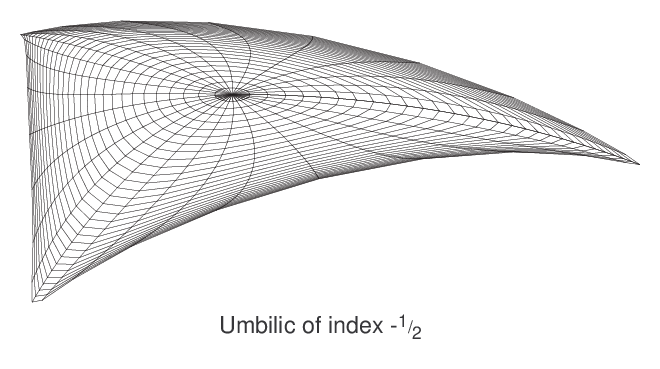}
  \end{center}
\end{wrapfigure}

Let us now construct the cross-cap. Consider the map $h:[R_0,1]\times S^1\rightarrow TS^2$, which takes $\nu=Re^{i\theta}$ to
$(\xi,\eta)=((1-\nu\bar{\nu})\nu,\bar{\nu}^2)$. 
This an embedded cross-cap since the image of $\{1\}\times S^1$ is a circle with antipodal points identified.

These two surfaces can be joined along a circle by fixing a constant Gauss radius $R_0$ such that $3^{-{\scriptstyle{\frac{1}{2}}}}<R_0<1$, and 
defining the surface $\tilde{\Sigma}=\Sigma\#{\mathbb R}P^2$ by
\begin{align}
\xi=&(1-\nu\bar{\nu})\nu \qquad\qquad {\mbox{ for}}\qquad 3^{-{\scriptstyle{\frac{1}{2}}}}<|\nu|\leq 1,\nonumber\\
& \nonumber\\
\eta=&\left\{
\begin{array}{lcl}
(1+(1-\nu\bar{\nu})^2\nu\bar{\nu})^2(1-\nu\bar{\nu})^2\bar{\nu}^2  &{\mbox{ for}} & 3^{-{\scriptstyle{\frac{1}{2}}}}<|\nu|\leq R_0, \\
 & \\
(a+b(1-\nu\bar{\nu})+c(1-\nu\bar{\nu})^2)\bar{\nu}^2&{\mbox{ for}} & R_0\leq|\nu|\leq 1.
\end{array}
\right.\nonumber
\end{align} 
This surface is $C^2$-smooth for the following choice of constants
\begin{align}
a=&-(1-R_0^2)^4(5+2R_0^2-46R_0^4+54R_0^6-21R_0^8),\nonumber\\
&\nonumber\\
b=&-2(1-R_0^2)^3(6-51R_0^4+61R_0^6-24R_0^8),\nonumber\\
&\nonumber\\
c=&-6+18R_0^2+42R_0^4-180R_0^6+225R_0^8-126R_0^{10}+28R_0^{12}.\nonumber
\end{align}
Moreover, it is easily seen to be totally real and may be smoothed to $C^{2+\alpha}$.

In ${\mathbb R}^3$, this surface can be visualized as a 2-parameter family of oriented lines, or a 1-parameter family of ruled surfaces which doubly covers a cylinder.

\begin{wrapfigure}{l}{0.5\textwidth}
  \begin{center}
    \includegraphics[width=0.47\textwidth]{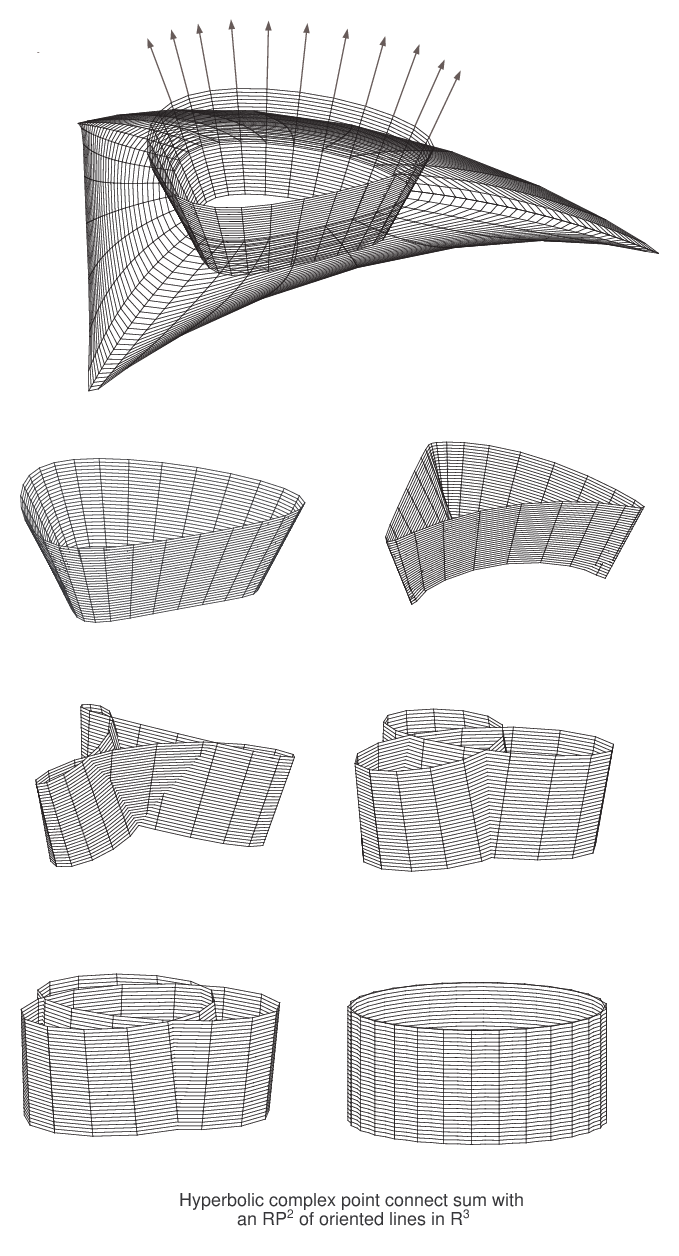}
  \end{center}
\end{wrapfigure}

We illustrate this with a sequence of the ruled surfaces corresponding to circles of constant Gauss radius $|\xi|$ on the cross-cap.
 
Removing the hyperbolic complex point means that the sum of the indices of the complex points $\gamma_j$ on $\Sigma$ has increased by one. This can be seen in ${\mathbb R}^3$ as follows. 

Recall Lai's formula for this sum:
$\sum_{j}I(\gamma_j)=\chi(T\Sigma)+\chi(N\Sigma)$.
Adding a cross-cap reduces $\chi(T\Sigma)$ by one and so a totally real blow-up must increase $\chi(N\Sigma)$ by two. 

To see that this is indeed the case, note that $\chi(N\Sigma)$ can be identified with the number of intersection points of $\Sigma$ with a small perturbation of $\Sigma$. Perturbing the cross-cap by a translation in ${\mathbb R}^3$ orthogonal to the central cylinder (addition of a certain quadratic holomorphic section in $TS^2$), the ruling of the cylinder and its perturbation coincide at exactly two lines. The rest of the cross-cap does not intersect its perturbation and so, replacing the disc by the cross-cap increases $\chi(N\Sigma)$ by two, as claimed.

\vspace{0.1in}

\section{Proof of the Theorem 3}\label{s:3}

\subsection{Reformulation}
Suppose, for the sake of contradiction, that there exists an open $C^{3+\alpha}$ convex surface containing an isolated umbilic point of index $i=2+k/2$, for $k\geq 0$. Close up the surface to form a closed convex surface $S$ by joining it to an umbilic-free hemisphere. Here hemisphere refers to the range of the Gauss map of the surface $\pi_S:S\rightarrow S^2$ and umbilic-free hemispheres are simple to construct. 

This closed surface $S$ will have umbilic points in the annulus, which are generically isolated and whose indices sum to $-k/2$. In fact, these generic umbilic points will be those classified by Darboux in 1896 \cite{Dar}, what we term {\it elliptic} umbilic points of index $1/2$ and {\it hyperbolic} umbilic points of index $-1/2$ (in-line with the complex terminology).

The reformulation outlined in Section \ref{s:1} means that the set of oriented normal lines to this surface form a global Lagrangian section $\Sigma$ in $TS^2$ with one complex point of index $I=4+k$, and isolated elliptic and hyperbolic complex points whose indices sum to $-k$.

\subsection{Cancellation of Complex Points}  

By a small deformation of $\Sigma$ we cancel the elliptic complex points so that $k$ hyperbolic complex points remain. Such cancellation can be carried 
out in general for complex points on real surfaces in a complex surface (see for example Corollary 9.5.2 in \cite{forst}). However, in our case we 
want to ensure that the deformation remains Lagrangian. This follows from:

\begin{Lem}
Let $\Sigma$ be a smooth simply connected Lagrangian section with non-empty totally real boundary. Suppose that the sum of the complex indices 
on $\Sigma$ is zero.  Then there exists a smooth Lagrangian section $\Sigma'$ with the following properties:
\begin{enumerate}
\item $\Sigma=\Sigma'$ on $\Sigma-K$ for $K$ a compact set in the interior of $\Sigma$,
\item $\Sigma'$ is totally real.
\end{enumerate}
\end{Lem}
\begin{proof}
Let $U=\pi(\Sigma)$ and pull back the jet-bundles over $\Sigma$ to ones over $U$. 

In particular, the support function $r:U\rightarrow {\mathbb R}$ gives a section in the 2-jet bundle of 
${\mathbb R}\oplus Hom({\mathbb R},{\mathbb R})\oplus Hom({\mathbb R}^2,{\mathbb R}^2)$ over $U$. This section is locally given by the map
\[
\xi\mapsto (r,r_1,r_2,r_{11},r_{12},r_{21},r_{22}),
\]
where we set $\xi=x^1+ix^2$ and a subscript denotes differentiation. Sections that arise in such a manner are called {\it holonomic} - see for example 
\cite{spring}. A non-holonomic section is of the form 
\[
\xi\mapsto (r,X_1,X_2,Y_1,Y_2,Y_3,Y_4),
\]
for functions $\{X_j\},\{Y_j\}$ on $U$. The boundary is totally real and so $|\partial\partial r|^2\neq0$ on $\partial U$. Note that
\[
|\partial\partial r|^2=(r_{11}-r_{22})^2+(r_{12}+r_{21})^2.
\]
Since the index sum is zero, the winding number of $\partial\partial r$ around the boundary is zero. Thus there exists a non-holonomic section of this bundle such that $(Y_1-Y_4)^2+(Y_2+Y_3)^2\neq0$ on $U$ which agrees with the original section in a neighbourhood of the boundary.

Consider the differential relation $(Y_1-Y_4)^2+(Y_2+Y_3)^2\neq0$. The convex hull of the complement of this relation is the full fibre (i.e. the relation is ample) and so by the relative h-principle, there exists a holonomic section $\Sigma'$ which is totally real and agrees with the original section on the boundary.

The surface $\Sigma'$ has the properties (1) to (2) above.

\end{proof}

Using this Lemma we can cancel pairs of hyperbolic and elliptic complex points by a perturbation.  We thus arrive at a Lagrangian section $\Sigma'$ with an isolated complex point of index $4+k$ and $k$ hyperbolic complex points, which also contains a totally real Lagrangian hemisphere.

\subsection{Removal of Hyperbolic Points}

Now remove the complex points of index $-1$ from the surface by totally real blow-up, as described in Proposition \ref{p:trb}. That is, we modify $\Sigma'$ by removing a small disc containing the complex point of index $-1$, and attaching ${\mathbb R}P^2$ with a disc removed (a cross-cap).  

Removing each of the $k$ complex points yields a closed embedded surface $\Sigma_1=\Sigma'\#k{\mathbb R}P^2$ containing a single isolated complex point (of index $I=4+k$), which is Lagrangian outside of the $k$ copies of ${\mathbb R}P^2$.

In summary, the sum of the indices of the complex points $\gamma_j$ on $\Sigma_1$ is
\begin{align}
\sum_{j}I(\gamma_j)&=\chi(T\Sigma_1)+\chi(N\Sigma_1)\nonumber\\
&=\chi(T(\Sigma'\#k{\mathbb R}P^2))+\chi(N(\Sigma'\#k{\mathbb R}P^2))\nonumber\\
&=\chi(T\Sigma')-k+\chi(N\Sigma')+2k\nonumber\\
&=4+k,\nonumber
\end{align}
which is the index of the single complex point on $\Sigma_1$. 

\subsection{Global Argument}
A minor modification of the global arguments of \cite{gak19} can now be applied to $\Sigma_1$. We sketch the arguments.

\begin{Def}
Let
\[
{\mathcal S}_0=\{\Sigma\;|\; \Sigma\subset TS^2{\mbox{ embedded }} C^{2+\alpha}{\mbox{ surface containing the point (0,0)}}\;\} 
\]
and the space of sections
\[
\Gamma_0({\mathbb J}(T\Sigma))
=\left\{\;v\in\Gamma({\mathbb J}(T\Sigma))\;\;\left|\;\; v\in C^{1+\alpha},\;\;v(\gamma_0)=0,\;\right.\right\}.
\]
\end{Def}

Let $\Sigma_1=\Sigma'\#k{\mathbb R}P^2$ be the closed surface in $TS^2$ with a single complex point of index $4+k$ obtained as above, and by a translation and rotation suppose that the complex point lies at $(0,0)\in TS^2$.

The neighbourhood of $\Sigma_1$ in ${\mathcal S}_0$ can be modeled by the sections $\Gamma_0$:

\begin{Prop}
There exists $\epsilon>0$ and $\Phi:B_{\epsilon}(0)\subset\Gamma_0({\mathbb J}(T\Sigma_1))\rightarrow{\mathcal U}\subset{\mathcal S}_0$ so that ${\mathcal U}=\Phi(B_{\epsilon}(0))$ is a Banach manifold.
\end{Prop}
\begin{proof}
The proof follows along the lines of the arguments in \cite{gak19}.
\end{proof}

In fact, we are interested in Lagrangian variations of $\Sigma_1$.
\begin{Def}
Let
\[
\Gamma_0^{\mbox{lag}}({\mathbb J}(T\Sigma_1))
=\left\{\;v\in\Gamma({\mathbb J}(T\Sigma_1))\;\;\left|\;\; v\in C^{1+\alpha},\;\;v(\gamma_0)=0,\;\;d({\mathbb J}(v)\lrcorner\Omega)=0\;\right.\right\}.
\]
\end{Def}

\begin{Prop}
$\Gamma_0^{\mbox{lag}}({\mathbb J}(T\Sigma_1))\subset\Gamma_0({\mathbb J}(T\Sigma_1))$ is a Banach subspace.
\end{Prop}

\begin{Cor}
There exists $\epsilon>0$ and $\Phi:B_{\epsilon}^{\mbox{lag}}(0)\subset\Gamma_0^{\mbox{lag}}({\mathbb J}(T\Sigma_1))\rightarrow{\mathcal U}^{\mbox{lag}}\subset{\mathcal U}$ so that ${\mathcal U}^{\mbox{lag}}=\Phi(B_{\epsilon}^{\mbox{lag}}(0))$ is a Banach submanifold.
\end{Cor}

Now the argument proceeds as follows. We conclude that for an open dense set of ${\mathcal U}^{\mbox{lag}}$, the co-kernel of the Cauchy-Riemann operator is zero. However, by Theorem 68 of \cite{gak24} there exist holomorphic discs whose boundary lie on any totally real Lagrangian hemisphere, contradicting the surjectivity of the $\bar{\partial}$-operator. 

Thus no $C^{2+\alpha}$ Lagrangian surface with an isolated complex point of index $I=4+k$ for $k\geq 0$ exists. Equivalently, there does not exist a $C^{3+\alpha}$ convex surface containing an isolated umbilic point of index $i=2+k/2$ for $k\geq 0$.

\vspace{0.1in}
\section{Discussion}

Since the early 1920's, attempts to prove the global Carath\'edory Conjecture have sought to establish a local index bound, usually with the inclusion of the assumption of real analyticity on the surface (originally due to Hamburger \cite{Ham} and more recently by Ivanov \cite{Ivan}). Conflating the analytic bound with a much later conjectured bound of Loewner (related to \cite{Loew}), the umbilic index bound most often considered, even in the smooth case, is $i\leq 1$.

Theorem \ref{t:1} represents a reversal of these historical attempts. That is, the results in this paper establish a local index bound as a consequence of global restrictions on the $\bar{\partial}$-operator. Remarkably, this opens up a gap between what we claim is a sharp result in the smooth category (umbilic index less than two) with Hamburger's result in the real analytic category (umbilic index less than or equal to one). 

Thus, we are led to the possibility of isolated umbilic points of index $3/2$ on smooth, non-real analytic surfaces. The existence and implications of such {\it exotic} umbilic points are worth considering for the future. 

\vspace{0.1in}

\noindent{\bf Acknowledgements}: The first author would like to thank the Institut des Hautes \'{E}tudes Scientifiques for support during the development of this work and the second author would like to thank the Max Planck Institut f\"ur Mathematik Bonn for support during its completion.

\end{document}